# Defining and Estimating Intervention Effects for Groups that will Develop an Auxiliary Outcome

**Marshall M. Joffe, Dylan Small and Chi-Yuan Hsu**

*Abstract.* It has recently become popular to define treatment effects for subsets of the target population characterized by variables not observable at the time a treatment decision is made. Characterizing and estimating such treatment effects is tricky; the most popular but naive approach inappropriately adjusts for variables affected by treatment and so is biased. We consider several appropriate ways to formalize the effects: principal stratification, stratification on a single potential auxiliary variable, stratification on an observed auxiliary variable and stratification on expected levels of auxiliary variables. We then outline identifying assumptions for each type of estimand. We evaluate the utility of these estimands and estimation procedures for decision making and understanding causal processes, contrasting them with the concepts of direct and indirect effects. We motivate our development with examples from nephrology and cancer screening, and use simulated data and real data on cancer screening to illustrate the estimation methods.

*Marshall M. Joffe is Associate Professor, Department of Biostatistics and Epidemiology, University of Pennsylvania School of Medicine, Philadelphia, Pennsylvania 19104-6021, USA e-mail: mjoffe@mail.med.upenn.edu. Dylan Small is Assistant Professor, Department of Statistics, The Wharton School, University of Pennsylvania, Philadelphia, Pennsylvania 19104-6340, USA e-mail: dsmall@wharton.upenn.edu. Chi-Yuan Hsu is Assistant Professor in Residence, Division of Nephrology, University of California, San Francisco, San Francisco, California 94143, USA e-mail: hsuchi@medicine.ucsf.edu.*

## 1. INTRODUCTION

In the recent literature on causal inference, it has become popular to define treatment effects for subsets of the target population characterized by variables not observable at the time a treatment decision is made. The most popular framework for doing this is principal stratification (PS); this name was introduced in a unifying paper by Frangakis and Rubin (2002). The ideas have been applied to a broad range of problems, including censoring by death (Robins, 1986; Zhang and Rubin, 2003), noncompliance in randomized trials (Angrist, Imbens and Rubin, 1996; Baker and Lindeman, 1994), the estimation of the effects of vaccines on post-infection outcomes (Gilbert, Bosch and Hudgens, 2003) and surrogate outcomes in randomized trials (Frangakis and Rubin, 2002). As we shall see, PS is one of several possible ways to define these effects.

The reasons for interest in effects defined by post-treatment auxiliary variables are diverse. We consider two problems, one in nephrology and one in cancer screening, to provide motivation for interest in these various estimands and estimation procedures.

Nephrologists have been frustrated by the lack of a good means to lower the high rates of morbidity and mortality among patients after the onset of end-stage renal disease (ESRD), the point at which kidney disease has progressed sufficiently to require dialysis treatment. Suggested methods to reduce this morbidity and mortality, which have included better or more aggressive control of anemia and higher doses of dialysis, have shown disappointing results (Paniagua, Amato, Vonesh, Correa-Rotter, Ramos, Moran and Mujais, 2002; Eknoyan, Beck, Cheung, Daugirdas, Greene, Kusek, Allon, Bailey, Delmez, Depner, Dwyer, Levey, Levin, Milford, Ornt, Rocco, Schulman, Schwab, Teehan and Toto, 2002; Besarab, Bolton, Browne, Egrie, Nissenson, Okamoto, Schwab and Goodkin, 1998). Frustration with the inability of treatments or interventions given after the initiation of dialysis to affect the course of ESRD has led some to hypothesize that the period before ESRD develops may provide a window of opportunity to improve outcomes among ESRD patients. The hope is that interventions or treatments applied before ESRD develops may affect the clinical course after development of ESRD; that is, it is hypothesized that treatments given before ESRD affect outcomes after the development of ESRD.

This apparently simple hypothesis is surprisingly hard to formalize. The difficulty stems from the fact that not all subjects with advanced chronic renal insufficiency (CRI) progress to ESRD, and that the same interventions that may affect outcomes after the onset of ESRD may themselves help determine who will develop ESRD. After considering a naive approach and illustrating its difficulties (Section 2), this paper considers several ways to formalize the hypothesis; for each, we briefly consider approaches to estimation of relevant and causally meaningful parameters (Section 3).

A second motivating problem concerns evaluation of the efficacy of cancer screening. Successful methods for screening for cancer result in earlier diagnosis of cancer (or precancerous conditions); this early detection may lead to treatment of the cancer while the cancer is still curable and so to reduced mortality. Randomized trials have been used to evaluate cancer screening; often, people randomized to screen fail to comply with their assignment. It is expected that any benefit of assignment to screening would be restricted to women who are screened. It might further be surmised that any benefit of screening would be restricted to screened women diagnosed with breast cancer and that the benefit of screening would be restricted even further to women whose cancer was diagnosed as a result of the screen. For explanatory purposes, it is of interest to estimate the benefit of screening for women in these subgroups. Hypotheses here may be formulated in ways similar to those in the nephrology problem. In these data, the outcome is failure-time, a censored continuous variable.

In Section 4 we consider statistical inference. We consider ranges of assumptions that identify the various estimands, as well as methods of estimation. Where necessary, we concentrate on continuous outcomes, as some methods for inference are more straightforward here. In Section 5 we consider a simulation experiment to provide some comparison of inference under various approaches. In Section 6 we analyze data about cancer screening from the Health Insurance Plan (HIP) Study (Shapiro, Venet, Strax and Venet, 1988), considering various estimands of interest.

In Section 7 we evaluate and compare the various approaches in terms of their utility for decision making and explanation. In addition, we compare the approaches conditioning treatment effects on post-treatment auxiliary variables to the estimation of direct and indirect effects. The paper concludes with discussion of extensions of the estimands and methods to more complex settings.

## 2. MOTIVATING EXAMPLE: DATA AND SIMPLE ANALYSIS

### 2.1 Data and Potential Outcomes

We motivate our methodological development with a simple numerical example. The example is from a study in a cohort of subjects with chronic renal insufficiency, a condition in which subjects have diminished kidney function but do not yet require dialysis or transplant (Feldman, Appel, Chertow, Cifelli, Cizman, Daugirdas, Fink, Franklin-Becker, Go, Hamm, He, Hostetter, Hsu, Jamerson, Joffe, Kusek, Landis, Lash, Miller, Mohler, Muntner, Ojo, Rahman, Townsend and Wright, 2003); although the numbers are arbitrary, they are intended to represent, in simplified fashion, the problems and associations present in studying effects in this population.

Here, we wish to study the effect of aggressive treatment of hypertension on myocardial infarction (MI) among subjects who develop ESRD after the start of follow-up. We make several simplifying assumptions, which we do not necessarily expect to apply in real data:

1. each of our main study variables is binary and scalar;
2. subjects are randomly assigned to either receive or not receive aggressive management of hypertension;
3. no subject has ESRD at the start of follow-up; and
4. subsequent ESRD status is recorded before any MI occurs.

We use $A$ to refer to the treatment of interest [$A = 1$ (0) indicates the presence (absence) of aggressive treatment of hypertension], $S$ to denote a post-treatment auxiliary variable [$S = 1$ (0) indicates the presence (absence) of ESRD] and $Y$ to refer to the outcome of interest [$Y = 1$ (0) indicates the occurrence (absence of occurrence) of an MI before the end of follow-up].

We adopt the potential outcomes approach (Neyman, 1990; Rubin, 1974) to illustrate and define our causal estimands of interest. Let $Y^a$ denote the outcome that would be seen were a subject given treatment level $a$ and let $S^a$ denote the level of the auxiliary variable were a subject given treatment $a$. Causal effects are normally defined in terms of comparisons of the outcomes that would be seen in the same individuals or groups under different conditions; for example, as comparisons of $Y^a$ and $Y^{a'}$ for $a \neq a'$. For this illustration, we assume that aggressive treatment does not affect MI for any individual (i.e., that $Y^0 = Y^1$ for all subjects); however, aggressive treatment will prevent ESRD for some subjects but never cause it, and so $S^1 \leq S^0$.

Table 1 classifies the population according to whether they would develop ESRD if treated and if untreated, and considers the risk of failure in all strata based on the cross-classification of this auxiliary variable. Half of the population would not develop ESRD whether or not they were treated (i.e., $S^0 = S^1 = 0$); the risk of MI in this group is low (10%). In another 30% of the population, aggressive treatment would prevent ESRD (i.e., $S^1 = 0, S^0 = 1$); the risk of MI in this group is higher (20%). The risk is highest (30%) in the 20% of the population doomed to get ESRD regardless of treatment (i.e., $S^1 = 1, S^0 = 1$).

In any study, we cannot simultaneously observe what would happen to any individual under aggressive treatment ($Y^1, S^1$) and in its absence ($Y^0, S^0$), and so the joint distribution of the variables represented in Table 1 is not estimable. Table 2 shows what would be observed in a randomized trial in which half of the subjects receive aggressive treatment and half do not. Because 30% of subjects in the cohort would develop ESRD only if not treated aggressively (Table 1, row 2), the proportion of subjects treated aggressively who develop ESRD (20% = 100/500) is much lower than the proportion of subjects not treated aggressively who develop ESRD (50% = 250/500). Among untreated subjects who develop ESRD, 24% develop an MI, whereas 30% of treated subjects who develop ESRD also develop an MI. Similarly, 10% of untreated subjects who do not develop ESRD later develop an MI, whereas 13.75% of treated subjects who do not develop ESRD later develop an MI. All of this may be derived (in expectation) from Table 1.

## 2.2 A Naive Approach

To examine the effect of aggressive treatment on MI among subjects who develop ESRD, some would compare the probability of MI among aggressively treated subjects who develop ESRD (0.3) with the probability of MI among untreated subjects who develop ESRD (0.24). A naive interpretation is that the difference between these probabilities (0.06) represents the effect of treatment for subjects who will develop ESRD. This interpretation is not correct, since aggressive treatment actually has no effect on MI for any subject. The naive comparison $pr(Y=1|A=1, S=1) - pr(Y=1|A=0, S=1)$ diverges from the true individual causal effects because membership in the groups being compared depends on treatment. Aggressive treatment reduces the number of subjects with ESRD. Subjects who would develop ESRD even if treated aggressively comprise the sickest subgroup in the study. A higher proportion of them than any other subgroup would have developed MI whether or not they had received aggressive treatment, and comparing them with untreated subjects who develop ESRD is comparing them to a combination of subjects from the same subgroup and subjects who would develop ESRD only if not treated aggressively. Thus, conditioning on ESRD, a post-treatment variable, leads to inappropriate estimates of overall treatment effect (Robins, Blevins, Ritter and Wulfsohn, 1992; Rosenbaum, 1984).

TABLE 1
*Probability of MI if treated aggressively or not, by ESRD status if treated aggressively or not*

| ESRD status | | Probability of MI | | |
| if untreated ($S^0$) | if treated ($S^1$) | $pr(Y^0 = 1 \mid S^0, S^1)$ | $pr(Y^1 = 1 \mid S^0, S^1)$ | $N$ |
|---|---|---|---|---|
| 0 | 0 | 0.1 | 0.1 | 500 |
| 1 | 0 | 0.2 | 0.2 | 300 |
| 1 | 1 | 0.3 | 0.3 | 200 |

TABLE 2
*Probability of MI in randomized trial, by treatment arm and observed ESRD status*

| Treatment ($A$) | ESRD ($S$) | $N$ | Probability of MI ($pr(Y=1|S,A)$) |
|---|---|---|---|
| 0 | 0 | 250 | 0.1 |
|   | 1 | 250 | 0.24 |
| 1 | 0 | 400 | 0.1375 |
|   | 1 | 100 | 0.3 |

A complementary approach to view the bias resulting from conditioning on ESRD involves directed acyclic graphs (DAGs) (Pearl, 1995, 2000). Figure 1 shows the relations between the observed variables $A, S$ and $Y$, and unobserved variable(s) $U$. The arrows from $U$ to both $S$ and $Y$ indicate that there are some unmeasured common causes of both. The arrow from $A$ to $S$ indicates that $A$ influences $S$, whereas the absence of any directed path (i.e., a series of directed arrows) from exposure to outcome indicates that $A$ has no effect on $Y$. $S$ is known as a collider (i.e., there are arrows which converge on $S$). It is well known that conditioning on colliders induces associations between the parents of the collider (i.e., between $A$ and $U$). Because $U$ influences $Y$, the conditional association between $A$ and $U$ given $S$ propagates to a conditional association between $A$ and $Y$ given $S$. Structurally similar problems involving selection bias in epidemiology are discussed elsewhere (Greenland, 2003; Hernán, Hernández-Díaz and Robins, 2004).

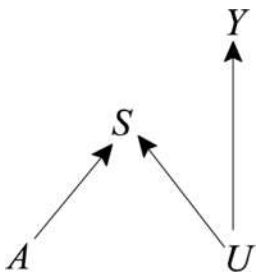

FIG. 1. *A directed acyclic graph representing the relations among the variables in the example of Section 2.*

We consider briefly an appropriate causal interpretation of the naive comparison of observable conditional distributions $pr(Y=1|A=1, S=1) - pr(Y=1|A=0, S=1)$. In a randomized trial, the observed distribution $pr(Y|A=a, S=s)$ equals the conditional distribution $pr(Y^a|S^a=s)$ of the potential outcome that would be seen among subjects who would have a common value of the auxiliary variable $S$ if they receive level $a$ of treatment. Here, this is the probability of MI that would be seen in subjects developing ESRD were all subjects to be treated aggressively (for $a=1$) or not receive aggressive treatment (for $a=0$). A comparison of $pr(Y^a|S^a=s)$ for different values of treatment $a$ reflects the impact of treatment on the conditional distribution of the potential outcome given the auxiliary outcome. In nonrandomized studies, the observed conditional distributions $pr(Y|A=a, S=s)$ will not in general equal the potential conditional distributions $pr(Y^a|S^a=s)$. So long as information is collected on all subjects in a cohort, this conditional distribution is identified under the commonly used assumption of strongly ignorable treatment assignment (Rosenbaum and Rubin, 1983),

(1)     $pr(A=a|\underline{S}, \underline{Y}, X) = pr(A=a|X) > 0,$

where $\underline{S} \equiv \{S^a\}$ and $\underline{Y} \equiv \{Y^a\}$ denote the vectors of potential auxiliary and main outcomes, respectively. Strong ignorability identifies $pr(Y^a|X, S^a = s)$ as $pr(Y|X, A=a, S=s)$. The potential conditional distributions may be useful in decision making (Section 7.1).

Sometimes, the effects of an intervention on a (conditional) distribution will be all that is identifiable from one's data under plausible assumptions. Here, however, there are other measures of effect that more closely relate to the scientific questions of interest, regarding the effect that aggressive treatment has for subjects who develop ESRD. We consider these in the following sections.

## 3. DEFINITION OF EFFECTS FOR A COMMON SET OF INDIVIDUALS

We consider several ways of characterizing or defining effects for subjects who will develop (or are likely to develop) an auxiliary outcome. In each approach, we consider first the group or subgroup for whom effects are defined; for each approach, the definition of effects for this group follows in a straightforward way from the potential outcomes approach. We relate the effects to each other through a common probability model in Section 3.7.

### 3.1 Principal Stratification

Principal stratification (PS) (Frangakis and Rubin, 2002) is a method proposed recently by several authors for defining certain types of causal effects. In this approach, effects are characterized within strata defined by the vector of potential auxiliary outcomes (here $S^0$ and $S^1$). In our example, we may be interested in the effect of aggressive treatment for subjects who would develop ESRD whether or not they are treated aggressively (i.e., $S^0 = S^1 = 1$). In particular, we concentrate here on comparing the proportion in this subgroup that would have an MI if treated aggressively $[pr(Y^1 = 1|S^0 = S^1 = 1)]$ and the proportion in the same subgroup who would have an MI if not treated aggressively $[pr(Y^0 = 1|S^0 = S^1 = 1)]$. In general, the approach compares the expected value or distribution of potential outcome $Y^a$ for different levels of treatment ($a = 0, 1$) in strata defined by the levels the auxiliary variable would take under both levels of treatment (i.e., strata are defined by $S^0$ and $S^1$ jointly). In our data, this stratum is shown in the third data row of Table 1; aggressive treatment has no effect on MI in this group (as in all other principal strata).

In general, the principal strata are not fully identified from the data, because one cannot simultaneously observe both potential auxiliary outcomes $S^0$ and $S^1$. In the data in Table 1, aggressive treatment sometimes prevents but never causes ESRD (i.e., $S^0 \geq S^1$); this monotonicity (Angrist, Imbens and Rubin, 1996) is not completely plausible in the nephrology example; Section 8.1 discusses this at more length. Under monotonicity, the principal stratum of some subjects is identifiable: aggressively treated subjects who develop ESRD would have done so even had they not been treated (i.e., $A = 1$, $S = S^1 = 1$ implies $S^0 = S^1 = 1$), and untreated subjects who do not develop ESRD would not have done so even had they been treated (i.e., $A = 0$, $S = S^0 = 0$ implies $S^1 = S^0 = 0$).

### 3.2 Single Potential Stratification

One can define the effect of a treatment for a subgroup defined by a single potential outcome. For example, one may be interested in the effect of aggressive treatment for people who would develop ESRD if they were treated aggressively [a comparison of $pr(Y^1 = 1|S^1 = 1)$ and $pr(Y^0 = 1|S^1 = 1)$; Table 1, row 3], or in the effect of aggressive treatment on people who would develop ESRD if they were not treated aggressively [a comparison of $pr(Y^1 = 1|S^0 = 1)$ and $pr(Y^0 = 1|S^0 = 1)$; the last two rows in Table 1]. Such stratification may be viewed as a coarser form of PS.

Membership in this stratum, defined by a single auxiliary variable ($S^0$ or $S^1$), is only partially observed. Whether a person is in this single auxiliary stratum defined by $S^a$ is known if a person receives treatment level $a$. As above, this complicates statistical inference for effects defined for such groups. One may use approaches such as those used for PS to estimate effects defined in this fashion; because of the connection with observed auxiliary stratification, one can also use methods described for the next sort of stratification we discuss.

### 3.3 Observed Auxiliary Stratification

For this approach, we define effects for groups defined by observed auxiliary variables. For example, we might consider the effect of aggressive treatment for subjects who received nonaggressive treatment and developed ESRD; that is, we compare $pr(Y^1 = 1|S = 1, A = 0)$ with $pr(Y^0 = 1|S = 1, A = 0)$. The subjects for whom this effect is defined are the 250 subjects in the second row of Table 2 (who comprise the 50% of the subjects in the second and third rows of Table 1 who are untreated). Alternatively, we might be interested in the effect of aggressive treatment for subjects who received aggressive treatment and developed ESRD; that is, we compare $pr(Y^1 = 1|S = 1, A = 1)$ with $pr(Y^0 = 1|S = 1, A = 1)$ (here, the effects are defined for the last row of Table 2, or the 50% of the last row of Table 1 who are treated). We have called the latter comparison the *realized effect* of treatment (Joffe, 2001), because it represents the effect treatment actually had within a subgroup (which, in this instance, is defined by post-treatment variables).

In a randomized trial, observed auxiliary stratification is equivalent to single potential auxiliary stratification. To see this, note that $pr(Y^a = 1|S = 1, A = a) = pr(Y^{a'} = 1|S^a = 1, A = a) = pr(Y^{a'} = 1|S^a = 1)$; the last step follows because of randomization. Thus, for example, in a randomized trial in our ESRD example, aggressively treated subjects who develop ESRD are comparable to the set of subjects who, if they had received aggressive treatment, would have developed ESRD. In observational studies, under ignorable treatment assignment (Rosenbaum and Rubin, 1983), the above holds conditional on covariates.

Here, the subgroup for whom the effect is defined is fully observable; however, the effects are defined for groups that are not identified at time of treatment decision. Thus, like PS and single potential stratification, this approach cannot be used directly to predict, at the time of a treatment decision, the effect of that decision. Further, there is an explanatory flavor to the analysis and model: these effects, which have already happened to defined subgroups, can be used to explain differences between randomized groups.

### 3.4 Expected Auxiliary Stratification

Unlike previous approaches, one can define effects in a way that uses information on auxiliary variables to define subgroups identifiable at the time treatment decisions are made. Here, we define effects for a group of subjects who are likely to develop ESRD if given a particular treatment. Let $\mu^a(X) \equiv E(S^a|X)$ denote the probability, given baseline covariates $X$, of developing ESRD if one were to receive treatment level $a$; it is the expected value of the potential auxiliary variable. $\mu^a(X)$ has been called a "principal score" (Hill, Waldfogel and Brooks-Gunn, 2002). We then define effects as a comparison of potential outcomes (MI) for subjects with the same expected auxiliary (ESRD), for example, $E\{Y^1|\mu^1(X)\} - E\{Y^0|\mu^1(X)\}$. Alternatively, we can define effects for broader strata based on $\mu^a(X)$, for example, $E\{Y^1|\mu^1(X) > 0.8\} - E\{Y^0|\mu^1(X) > 0.8\}$, the effect of aggressive treatment for subjects with at least an 80% chance of developing ESRD if treated aggressively.

### 3.5 Expected Multiple Auxiliary Stratification

The approach of the last section may be extended to condition effects on multiple expected auxiliary variables. Thus, one might be interested in the effect of treatment for subjects with an 80% risk of developing ESRD if treated aggressively and a 90% risk if not treated aggressively; that is, we derive effects for subjects based in groups determined by both $\mu^0(X)$ and $\mu^1(X)$. This approach has the flavor of PS; unlike PS, the subgroups for whom effects are defined are fully identified in the data, based solely on pretreatment information. As in the previous section, conventional statistical methods apply.

### 3.6 Conventional Approach

A final alternative is to estimate effects in subgroups based solely on pretreatment covariates, where group membership is not dependent on any risk score like $\mu^a(X)$. Conventional estimation approaches may apply; we may look for effect modification by baseline covariates directly, rather than indirectly through the expected auxiliary.

### 3.7 Probability Models for the Data

Following Rubin (1978), we consider a formal probability model for the joint distribution of the observable data and potential outcomes. We then consider approaches to parametrizing parts of this distribution (Section 4).

One general way to factor the joint density of the observable quantities and potential outcomes is

$$(2) \quad \begin{aligned} &f(X, \underline{S}, \underline{Y}, A) \\ &= f(X)f(\underline{S}|X)f(\underline{Y}|X, \underline{S})f(A|X, \underline{S}, \underline{Y}). \end{aligned}$$

This factorization is akin to selection models in common use in longitudinal data analysis (Little, 1995). Strongly ignorable treatment assignment (Rosenbaum and Rubin, 1983) is often assumed.

PS estimands depend only on one part of the joint density: $f(\underline{Y}|X, \underline{S})$. Other causal estimands also involve the density of the auxiliary outcome (and sometimes of the exposure). Single potential stratification estimands involve $f(\underline{Y}|X, S^a)$. These estimands can be obtained by integrating out the other potential auxiliary outcomes; that is, $f(\underline{Y}|X, S^a) = \int_{\underline{S}^{\neg a}} f(\underline{Y}|X, \underline{S})f(\underline{S}|X) d\underline{S}^{\neg a}/f(S^a|X)$, where $\underline{S}^{\neg a}$ refers to the vector of unobserved potential auxiliary outcomes. Observed auxiliary estimands involve $f(\underline{Y}|X, S, A) = f(\underline{Y}, A|X, S)/f(A|X, S) = \int_{\underline{S}^{\neg a}} f(\underline{Y}|X, \underline{S})f(\underline{S}|X) \cdot$

$f(A|X, \underline{S}, \underline{Y}) \, d\underline{S}^{\neg a} / \int_{\underline{S}^{\neg a}} \int_{\underline{Y}} f(\underline{Y}|X, \underline{S}) f(\underline{S}|X) f(A|X,$
$\underline{S}, \underline{Y}) \, d\underline{S}^{\neg a} \, d\underline{Y}$; under ignorability, this simplifies to $f(\underline{Y}|X, S, A = a) = \int_{\underline{S}^{\neg a}} f(\underline{Y}|X, \underline{S}) f(\underline{S}|X) \, d\underline{s}^{\neg a} / f(S^a|X)$, which is the estimand of single potential stratification. Expected auxiliary stratification involves integrating out all the auxiliary potential outcomes; that is,

$$f\{\underline{Y}|E(S^a|X) = s\}$$
$$= \frac{\int_{X:E(S^a|X)=s} \int_{\underline{S}} f(\underline{Y}|X, \underline{S}) f(\underline{S}|X) f(X) \, d\underline{S} \, dX}{\int_{X:E(S^a|X)=s} f(X) \, dX}$$
$$= \frac{\int_{X:E(S^a|X)=s} f(\underline{Y}|X) f(X) \, dX}{\int_{X:E(S^a|X)=s} f(X) \, dX};$$

in this formulation, one can ignore the density of the auxiliary outcome except as it relates to identifying the subset over whom to average. Similarly, in stratifying on observed variables, one can integrate out and ignore the auxiliary outcome: $f(\underline{Y}|X) = \int_{\underline{S}} f(\underline{Y}|X, \underline{S}) f(\underline{S}|X) \, d\underline{X}$. Typically, inference concentrates on comparisons of the marginal densities $f(Y^a|X, \cdot)$ of the potential outcomes for different treatment levels $a$, rather than the joint density $f(\underline{Y}|X, \cdot)$ of the potential outcomes $\underline{Y}$.

## 4. CONTINUOUS OUTCOMES: MODELS AND ESTIMATION

This section considers statistical inference about the various estimands outlined in the previous section. Some estimation methods in this setting are more straightforward with continuous outcome data, and so we discuss models and estimation methods more fully for such data. We apply these methods both to simulated data (Section 5) and the HIP data (Section 6).

We now discuss estimation of the various causal quantities defined in the previous sections, concentrating on identification of causal contrasts. We reverse the order of discussion, beginning with stratification on observed variables and ending with PS, as the nature of the latent structure and identifying assumptions becomes increasingly complex. The greater the degree of latent structure, the more assumptions are needed to estimate parameters in the model. The greater degree of assumptions involved is justified if they lead to better decision making, more explanatory power or greater generalizability, issues we take up in Section 7.

### 4.1 Conventional Approaches

For the conventional approach, identification of the marginal distributions $f(Y^a|X)$ may be based on the assumption of ignorable treatment assignment (Rosenbaum and Rubin, 1983), $pr(A = a|X, \underline{Y}) = pr(A = a|X)$, with $0 < pr(A = a|X) < 1$ for all $X, a$. Under this assumption, the observed density $f(Y|X, A = a)$ equals the density of the potential outcome $f(Y^a|X)$. Thus, one can use standard approaches (e.g., regression of $Y$ on $X$ and $A$) to estimate the effect of $A$ on $Y$.

### 4.2 Stratification on Expected Auxiliaries

For stratification on expected auxiliary variables (Hill, Waldfogel and Brooks-Gunn, 2002), identification and estimation are somewhat more complicated. The simplest method of estimation involves two steps: estimating the expected auxiliary $\mu^a(X)$, and estimation of effects by level of this expected auxiliary. One can estimate the expected auxiliary under ignorable treatment assignment [here $f(A|X, \underline{S}) = f(A|X)$] (Rosenbaum and Rubin, 1983) using standard methods. For example, one can regress $S$ on $X$ for subjects with $A = a$, then compute the expected value as $\hat{\mu}^a(X)$ based on the estimated regression coefficients. Alternatively, one can regress $S$ on $X$ and $A$, then estimate $\hat{\mu}^a(X)$ as $E(S|X, A = a)$, using, for each subject, his or her observed $X$ and the desired treatment level $a$ for all subjects; here, one may choose to include interactions between $A$ and $X$ as appropriate.

Under ignorable treatment assignment for the outcome $Y$, one can again use standard methods to estimate the effect of treatment for a group classified by $\hat{\mu}^a(X)$. For example, one can fit a regression

$$(3) \quad \begin{aligned} &E\{Y|\hat{\mu}^a(X), A\} \\ &= \beta_0 + \hat{\mu}^a(X)\beta_{\hat{\mu}} + A\beta_A + \hat{\mu}^a(X)A\beta_{\hat{\mu}A}; \end{aligned}$$

here, the effect of treatment $E\{Y^1|\hat{\mu}^a(X)\} - E\{Y^0|\hat{\mu}^a(X)\}$ for subjects with expected auxiliary $\hat{\mu}^a(X)$ is $\beta_A + \hat{\mu}^a(X)\beta_{\hat{\mu}a}$. These plug-in type estimates are, in general, consistent for the true parameters $\beta_A$ and $\beta_{\mu a}$ in the corresponding regression on the true scores $\mu^a(X)$ but will typically be biased in small samples; this bias results from the fact that the estimate $\hat{\mu}^a$ is a mismeasured version on true expected auxiliary $\mu^a$. It is of some interest to develop unbiased estimators of these effects. Although we are unaware of any such work in this setting, such work has been done with other generated regressors (Pagan, 1984).

It is tempting to extrapolate the regression effect and interpret $\beta_A + \beta_{\mu a}$ as the effect for a subgroup in which all subjects will develop ESRD if treated aggressively, or as the effect for the group of subjects who will develop ESRD if treated aggressively. These interpretations are flawed. For the former, there may be no subgroups identifiable on the basis of pretreatment covariates in which all subjects will develop ESRD, and so the parameter is not meaningful. For the latter, additional assumptions may be required for this interpretation to hold, as has been discussed in a related context (Joffe, Ten Have and Brensinger, 2003).

### 4.3 Stratification on Observed Auxiliaries

We are unaware of previous work on estimation of parameters in models stratifying on observed auxiliary variables. Whereas the previous estimands allowed nonparametric identification under ignorability assumptions, estimation of these parameters will require additional assumptions. We sketch one approach to estimation below; because of the relation between observed auxiliary stratification and single potential stratification and PS (Section 3.3), estimation may also be based on estimating principal-strata specific effects as described in the next subsection, then marginalizing over the unobserved auxiliary outcomes $\underline{S}^{\neg a}$.

We consider estimation based on the following idea, similar to G-estimation in structural nested models (Robins, 1992; Robins et al., 1992; Robins, 1994). Estimation will require ignorability assumptions, as above. Suppose our outcome $Y$ is continuous; $Y$ could be the logarithm of a failure-time (e.g., time of breast cancer mortality). We first propose a model for the effect of treatment, for example,

$$
\begin{aligned}
(4) \quad & E(Y^0|X, A, S) \\
& = E(Y|X, A, S) - A(1-S)\Psi_0 - AS\Psi_1.
\end{aligned}
$$

Here, the realized effect of treatment for the subgroup of subjects who received treatment and developed ESRD is $\Psi_1$, and the realized effect for the subgroup which received treatment but did not develop ESRD is $\Psi_0$. We have assumed that these effects do not vary with covariates $X$.

Under ignorability, $Y^0$ is independent of $A$ given $X$. Let $Y^0(\Psi) = Y - A(1-S)\Psi_0 - AS\Psi_1$; $Y^0(\Psi)$ may be computed from observed quantities and putative values for the causal parameters $\Psi$. Based on the nonidentifiable assumption that treatment effects are the same for all subjects with common values of $A$ and $S$, $Y^0(\Psi)$ may be viewed heuristically as the potential outcome $Y^0$ if causal theories represented by $\{\Psi_0, \Psi_1\}$ are true. If the putative value of the causal parameter $\Psi$ is true, $Y^0(\Psi)$ will be independent of $A$ given $X$. Estimation may be based on testing this independence for an assumed value of the causal parameter $\Psi$. Because $\Psi$ is a vector of dimension 2, estimation using scalar estimating equations will require either restriction of the unknown parameter $\Psi$ or a vector of estimating equations of the same dimension as $\Psi$.

In some cases, one might assume that either $\Psi_0$ or $\Psi_1$ is 0. In the HIP study, it might be reasonable to assume that screening affects breast cancer mortality only for screened subjects diagnosed with breast cancer, or, further, only for screened subjects whose cancer was detected due to the screen (i.e., for subjects for whom $S = 1, A = 1$), and so $\Psi_0 = 0$. For binary $A$,

$$
(5) \quad \sum_i (A - p)g\{Y^0(\Psi), X\} = 0
$$

provides valid estimating equations under the model assumptions and ignorability, where $p \equiv pr(A = 1|X)$ and $g(\cdot)$ is a known function of its arguments. The optimal function $g(\cdot)$ is a sometimes complex function of the joint density of the observables and potential outcomes (Joffe and Brensinger, 2003; Robins, 1992; Robins et al., 1992). Efficiency, but not consistency, depends on choosing this optimal function. Suppose that the "error" terms $\varepsilon = Y^0 - E(Y^0|X)$ are normal, independent, identically distributed random variables and that $S^1$ is unrelated to $Y^0$ [i.e., $f(Y^0|X, A, S^1) = f(Y^0|X)$]. The optimal function is then $g\{Y^0(\Psi), X\} = \varepsilon(\Psi)E(S|X, A = 1)$, where $\varepsilon(\Psi) = Y^0(\Psi) - E\{Y^0(\Psi)|X\}$ (Joffe and Brensinger, 2003); the use of $\varepsilon(\Psi)$ is similar to Rosenbaum's (2002) use of such residuals in randomization-based inference. The estimated probability of treatment may be substituted for the typically unknown true $p$. The asymptotic variance of the resulting estimator may be derived using a sandwich-type formula (Robins, 1992; Robins et al., 1992).

If one is unable to restrict the parameter $\Psi$ based on subject-matter considerations, one must use a vector of estimating equations. Here, we will require vector functions $g\{Y(\Psi), X\}$. Under the above normality and homoscedasticity assumptions and under the assumption that $f(Y^0|X, A, S^1) = f(Y^0|X)$, the optimal function is the vector function $g\{Y^0(\boldsymbol{\Psi}), X\} = \varepsilon(\boldsymbol{\Psi})\{1 - E(S|X, A = 1), E(S|X, A = 1)\}^T$. $g\{Y^0(\boldsymbol{\Psi}), X\}$ must have the same rank as the dimension of $\boldsymbol{\Psi}$; thus, if the covariate does not predict the auxiliary outcome among the treated, there will be no ability to estimate the vector $\boldsymbol{\Psi}$. Under our assumptions of equal effect across strata of $X$, the effect of $A$ on $Y$ in any stratum of $X$ is $\boldsymbol{\Psi}_0\{1 - E(S|X, A = 1)\} + \boldsymbol{\Psi}_1 E(S|X, A = 1)$, each term corresponding to part of the vector function $g\{Y^0(\boldsymbol{\Psi}), X\}$.

The approach taken above is semiparametric; that is, consistency of the estimators of the parameters of interest $\boldsymbol{\Psi}$ does not depend on parametric assumptions about the distribution of $\varepsilon$ or the association of $X$ or $\underline{S}$ with $Y^0$. The methods presented above are valid for structural distribution models, which map percentiles in the distribution of $Y^1$ to percentiles in the distribution of $Y^0$ in defined subgroups. In contrast, the weaker structural mean models only map the mean of $Y^1$ to the mean of $Y^0$ [as in the formulation in (4)]. Consequently, there are fewer valid choices of the function $g(\cdot)$ for use in estimation [i.e., the function $g(\cdot)$ must be linear in $Y^0(\boldsymbol{\Psi})$] (Robins, Rotnitzky and Scharfstein, 2000).

With binary outcomes, mean models are required. Structural mean models using the identity link do not restrict subgroup-specific means to the permissible range. With the usual logit link, semiparametric models may be formulated, but consistent semiparametric estimators are not generally available (Robins, Rotnitzky and Scharfstein, 2000; Robins and Rotnitzky, 2004). Extensions to binary outcomes whose consistency depends on parametric assumptions are available (Robins and Rotnitzky, 2004; Vansteelandt and Goetghebeur, 2003).

### 4.4 Principal Stratification

Because the principal strata of many subjects are not determined by the usual combination of data and assumptions, inference for PS estimands will be more dependent on assumptions and complicated. The observed density of the main and auxiliary outcomes may be expressed in terms of the unobserved potential auxiliary outcomes $\underline{S}^A$ as follows:

$$f(Y, S|A = a, X)$$
$$= f(Y|S, A = a, X) f(S|A = a, X)$$
$$= \sum_{\underline{s}^{\neg a}} f(Y|S, \underline{S}^{\neg a} = \underline{s}^{\neg a}, A = a, X)$$
(6)
$$\cdot f(S^a, \underline{S}^{\neg a} = \underline{s}^{\neg a}|A = a, X)$$

$$= \sum_{\underline{s}^{\neg a}} f(Y|S, \underline{S}^{\neg a} = \underline{s}^{\neg a}, A = a, X)$$
$$\cdot f(\underline{S}^{\neg a} = \underline{s}^{\neg a}|S, A = a, X)$$
$$\cdot f(S|A = a, X).$$

Under ignorability, the components of the right-hand side of (6) are equivalent to models for the potential main and auxiliary potential outcomes, $f(Y|S, \underline{S}^{\neg a}, A = a, X) = f(Y^a|\underline{S}, X)$ and $f(S, \underline{S}^{\neg a}|A = a, X) = f(\underline{S}|X)$. Thus, parametrizing the models for the observables can lead to a likelihood for the causal quantities of interest.

This likelihood is generally overparametrized. Suppose that both $S$ and $A$ are scalar and binary. Then, there are six unknown observed densities for each level of covariates $X$: two for $f(S|A = a, X)$ (one for each level of $A$), and four for $f(Y|S, A = a, X)$ (one for each level of $A$ and $S$). However, there are eleven densities for the causal parameters in (6): three for $f(\underline{S}|X)$ (there are four levels of $\underline{S}$, but the probabilities sum to 1), and eight for $f(Y^a|X, \underline{S})$ (two levels of $a$ × two levels of $S$ × two levels of $\underline{S}^{\neg A}$). Thus, identification of causal effects will require further restrictions on the parameters. Before proceeding to discuss such restrictions, we note that there are (at least) two approaches which do not require identification: a Bayesian approach, in which prior information is combined with the likelihood to produce a posterior (Imbens and Rubin, 1997), and approaches which derive bounds on causal effects (Balke and Pearl, 1997; Manski, 1990; Robins, 1989; Rubin, 2004; Zhang and Rubin, 2003; Cheng and Small, 2006).

There are several forms of restrictions that can be applied. These include restrictions on the joint distribution of the potential auxiliary outcomes $\underline{S}$ (i.e., on the principal strata), and restrictions on the marginal distributions of the potential main outcomes $Y^a$.

One type of restriction on the auxiliary outcome occurs when some level of the auxiliary outcome is impossible under some treatment level $a$. Consider first randomized trials with noncompliance, one area in which the ideas of PS have been applied frequently (Angrist, Imbens and Rubin, 1996; Baker and Lindeman, 1994). In one view of these studies (Imbens and Rubin, 1997), one may view randomization, the controlled factor, as the treatment $A$ whose effect one is interested in estimating, and the level of exposure to the experimental therapy $S$ as the auxiliary variable. One may be interested in the effects of randomization for different classes of subjects defined by the set of behaviors $\underline{S} = \{S^0, S^1\}$ they would follow under different treatment assignments. In some randomized trials, subjects assigned the control treatment $(A = 0)$ have no access to the experimental therapy $(S = 1)$; this may be common with investigational therapies not available outside of the trial. In this case, there are two rather than four principal strata $\underline{S}: \{S^0, S^1\} = (0,1)$ (compliers) and $\{S^0, S^1\} = (0,0)$ (never-takers). In the HIP Study, it may be reasonable (under the conditions of the study in the 1960s) to reason that subjects not assigned to screening $(A = 0)$ could not have received mammograms and so could not have been diagnosed as a result of mammography $(S = 1)$. In these cases, the principal stratum of the treated subjects is known; since $S^0 = 0$, knowing the observed outcome $S = S^1$ in a treated subject fully identifies $\underline{S}$. Further, under ignorability, the proportion of subjects in each principal stratum $pr(\underline{S}|X)$ is identified. Under ignorability, the densities $f(Y^1|X, S^0, S^1 = s)$ are identified as $f(Y|X, S = s, A = 1)$, but the densities $f(Y^0|X, \underline{S}, A = 0)$ are not identified, and so the stratum-specific causal effects are not identified without further assumptions.

A more general form of restriction on the potential auxiliary outcomes consists of monotonicity restrictions, which have been adduced in a variety of settings (Efron and Feldman, 1991; Gilbert, Bosch and Hudgens, 2003; Angrist, Imbens and Rubin, 1996). For any two ordered levels of treatment $a, a', a' > a$, strict monotonicity states that $S^{a'} \geq S^a$ for all subjects (or, alternatively, that $S^{a'} \leq S^a$). Monotonicity identifies the principal stratum of some subjects. For binary $S$, monotonicity and ignorability together permit identification of the proportion of subjects in each principal stratum. In the HIP study, if $S$ is cancer diagnosis, one might assume that, if a cancer were detected in an unscreened subject, it also would have been detected had the subject been screened (i.e., $S^1 > S^0$). This assumption may be incorrect; for example, a woman who has been screened and who was told that she had no cancer may become less suspicious of cancer later, and so screening might sometimes lead to a missed diagnosis of cancer.

None of these monotonicity assumptions, even in conjunction with ignorability, is sufficient to identify stratum-specific treatment effects. Such identification typically requires assumptions about the outcome distributions $f(Y^a|X, \underline{S})$. We outline several such assumptions.

One strict set of assumptions is that treatment has no effect on the outcome $Y$ for some subsets of the data defined by $\underline{S}$. In randomized trials, it is common (Mark and Robins, 1993a; Sommer and Zeger, 1991) to assume that randomization has no effect for those subjects for whom treatment has no effect on $S$ [i.e., $Y^a = Y^{a'}$ if $S^a = S^{a'}$, or $f(Y^a|X, \underline{S}) = f(Y^{a'}|X, \underline{S})$ if $S^a = S^{a'}$, a weaker assumption]; this assumption is known as an exclusion restriction in the econometrics literature (Angrist, Imbens and Rubin, 1996). In studies of breast cancer screening (e.g., HIP), one might assume that screening has no effect among subjects whose tumors were not detected through screening (i.e., for subjects with $S^1 = 0$, $S^0 = 0$ where $S = 1$ if a subject has a screen-detected tumor, 0 otherwise). In both of these examples, these restrictions, together with the monotonicity restrictions and ignorability, are sufficient to identify the causal effects [and, in fact, the marginal densities $f(Y^a|X, \underline{S})$] in the single principal stratum in which treatment effects are not assumed to be 0. To see this, note that the marginal density of the outcome may be written $f(Y^a|X) = \sum_{\underline{s}} pr(\underline{S} = \underline{s}) f(Y^a|X, \underline{S} = \underline{s})$; under our assumptions, $pr(\underline{S} = \underline{s})$ is known for all $\underline{s}$, and $f(Y^a|X, \underline{S} = \underline{s})$ is identified as either $f(Y^a|X, S, A = a)$ or $f(Y^a|X, S, A = 1 - a)$ in the principal strata in which treatment has no effect. Further, the marginal densities $f(Y^a|X)$ are identified from the data under ignorability. This leaves one unknown quantity $f(Y^a|X, \underline{S} = \underline{s})$ for the principal stratum in which the effect is not assumed to be zero; we then solve $f(Y^a|X) = \sum_{\underline{s}} pr(\underline{S} = \underline{s}) f(Y^a|X, \underline{S} = \underline{s})$ to find this unknown density.

When these assumptions are not reasonable, the principal stratum-specific effects are not identified without other assumptions. These assumptions can take several forms: assumptions about how different principal strata are from each other with respect to the potential outcomes, assumptions about the distribution of potential outcomes within principal strata, and assumptions about the associations of covariates $X$ with the potential outcomes within principal strata.

Assumptions about the differences between principal strata may involve specifying $E(Y^0|X, \underline{S} = \underline{s}) - E(Y^0|X, \underline{S} = \underline{s}')$ for some $\underline{s}, \underline{s}', \underline{s} \neq \underline{s}'$; such assumptions are also sometimes framed in terms of the degree of dependence of $S^{1-A}$ on $Y^0$ in models for

$pr(S^{1-A}|X,S,A,Y^0)$. Such assumptions allow one to avoid parametric assumptions about error distributions. It is fairly difficult to have precise quantitative knowledge about the degree of this dependency. One way to deal with such uncertainty is to perform a sensitivity analysis, in which one varies the sensitivity parameter(s) over a plausible range (Gilbert, Bosch and Hudgens, 2003); another way is to put a prior distribution on the dependency parameter, as done in various Bayesian analyses (Imbens and Rubin, 1997).

Assumptions about the distribution of outcomes $Y^a$ within a principal stratum are sometimes made in both Bayesian and frequentist inference. Most often, $f(Y^a|\underline{S})$ is assumed to be normal with unknown mean and variance (Imbens and Rubin, 1997); we will consider likelihood-based inference under this assumption. With these assumptions, identification results from the identifiability of mixtures of certain parametric families of distributions (Teicher, 1963). For example, consider a model with two principal strata, never-takers and compliers, for which no exclusion restriction is assumed but the densities $f(Y^a|\underline{S})$ are assumed to be normal. Identifiability of $f\{Y^0|\underline{S}=(0,1)\}$ and $f\{Y^0|\underline{S}=(0,0)\}$ results from the fact that $f(Y|A=0)$ is a mixture of its two component normal distributions, and the parameters of a mixture of two normals are identifiable (Teicher, 1963). In order to identify which mixture component represents $f\{Y^0|\underline{S}=(0,1)\}$ and which represents $f\{Y^0|\underline{S}=(0,0)\}$, it is required that the proportion of compliers not equal 0.5. Identification based on parametric assumptions about the distributions of outcomes $f(Y^a|\underline{S})$ may not be robust to changes in parametric assumptions, as with selection models (Copas and Li, 1997).

Finally, one can make assumptions about the associations of covariates $X$ with the potential outcomes within principal strata. For example, one might assume that principal stratum-specific effects are constant across covariates $X$ [e.g., $E(Y^1-Y^0|X=x,\underline{S})=E(Y^1-Y^0|X=x',\underline{S})$, for all $x,x'$]; this assumption parallels assumptions we made for structural nested models that effects conditional on observed auxiliaries do not vary with $X$. One might also assume that the association of the covariates with potential outcomes has a particular parametric form [e.g., $E(Y^a|X,\underline{S})=q(\underline{S})+X\beta$].

## 5. A SMALL SIMULATION EXPERIMENT

In this section, we report results of a small simulation study to examine the performance of certain PS and observed auxiliary stratification estimators. We consider a continuous outcome $Y$, one binary covariate $X$, and four settings (two sets of expected values for outcomes within the $A,X,\underline{S}$ strata and two distributions of outcomes within the $A,X,\underline{S}$ strata for each set of expected values). For the four settings, the expected values of $Y$ within the $A,X,\underline{S}$ strata are shown in Table 3 along with the parametric distribution within the $A,X,\underline{S}$ strata. We label the principal strata $\underline{S}$ by $(S^0=0,S^1=0)=$ immune (I), $(S^0=1,S^1=0)=$ treatment protective (TP), $(S^0=0,S^1=1)=$ treatment harmful (TH) and $(S^0=1,S^1=1)=$ doomed (D) (Greenland and Robins, 1986). For all settings, the standard deviation of $Y$ within each $A,X,\underline{S}$ stratum is 1; $E(Y|A,X=1,\underline{S})-E(Y|A,X=0,\underline{S})=0.5$ for all $A,\underline{S}$, so that the principal stratum-specific effects are constant across the covariate $X$; $pr(\underline{S}=\mathrm{I})=0.25$, $pr(\underline{S}=\mathrm{TP})=0.4$, $pr(\underline{S}=\mathrm{TH})=0.05$ and $pr(\underline{S}=\mathrm{D})=0.3$; and $pr(X=1|\underline{S}=\mathrm{I})=0.5$, $pr(X=1|\underline{S}=\mathrm{TP})=0.75$, $pr(X=1|\underline{S}=\mathrm{TH})=0.25$ and $pr(X=1|\underline{S}=\mathrm{D})=0.5$. The sample size is 5000 for each setting and the probability of being randomly assigned treatment $(A=1)$ is 0.5.

We consider two PS estimators and one observed auxiliary stratification estimator. The PS estimators assume normal outcomes within each $A,X,\underline{S}$ stratum. This assumption is satisfied for settings I(A) and II(A) and violated for settings I(B) and II(B). The PS estimators furthermore make the correct assumption that $E(Y|A,X=1,\underline{S})-E(Y|A,X=0,\underline{S})$ is the same for all $A,\underline{S}$. The first PS estimator we consider does not put any constraints on $E(Y|A,X,\underline{S})$. The second estimator constrains the average effect of treatment in the immune principal stratum given a value of $X[E(Y|A=1,X,\underline{S}=\mathrm{I})-E(Y|A=0,X,\underline{S}=\mathrm{I})]$ to be equal to the average effect of treatment in the treatment protective principal stratum given the same value of $X$ and constrains the average effect of treatment in the treatment harmful principal stratum given a value of $X$ to be equal to the average effect of treatment in the doomed principal stratum given the same value of $X$. The constraint made by the second PS estimator is correct for settings I(A) and I(B) but is incorrect for settings II(A) and II(B). We used the EM algorithm to implement the PS estimators and used the true parameters as the starting values. The observed auxiliary stratification estimator is based on model (4) and uses the function $g\{Y^0(\Psi),X\}=\varepsilon(\Psi)\{1-E(S|X,A=1),E(S|X,A=$

$1)\}^T$ discussed in Section 4.3. Model (4) is correct for settings I(A) and I(B) but is incorrect for settings II(A) and II(B); in settings II(A) and II(B), $E(Y|X,A,S) - E(Y^0|X,A,S)$ varies with $X$. For the observed auxiliary stratification estimator, we used the computational procedure for weighted G-estimation described in Ten Have, Elliott, Joffe, Zanutto and Datto (2004) using the true values as the starting values. R code for the simulations and analysis and an example dataset are available at *www.cceb.upenn.edu/faculty/?id=157*.

Table 4 displays the means and standard deviations of the estimators for 500 simulations of each setting. In setting I(A), the assumptions made by all three estimators are correct. All the estimators provide approximately unbiased estimates of their corresponding estimands. The PS II estimator that constrains certain average effects of treatment within principal stratum to be equal is the most efficient estimator. The PS II estimator and the observed auxiliary stratification estimator make similar assumptions, but the additional parametric assumptions made by the PS II estimator make it considerably more efficient. In setting I(B), the assumptions made by the observed auxiliary stratification estimator continue to hold but the parametric assumptions made by the PS estimators are false. The observed auxiliary stratification estimator performs similarly in setting I(B) as it did in setting I(A)—it is approximately unbiased and has a similar variance. In contrast, the PS estimators exhibit considerable bias for some of their estimands in setting I(B). The PS I estimator has a particularly large bias for the average effects of treatment in the immune and doomed principal strata. In settings II(A) and II(B), the assumption made by the observed auxiliary stratification estimator that treatment effects conditional on observed auxiliaries do not vary with the covariate $X$ is false. In the "Truth" column for $\Psi_0$ and $\Psi_1$, we list the average realized effects of treatment for subjects with $S^1 = 0$ and $S^1 = 1$, respectively. The observed auxiliary stratification estimator of $\Psi_0$ for settings II(A) and II(B) does not show much bias but the estimator of $\Psi_1$ shows considerable bias. For the PS estimators, the assumptions made by PS estimator I are true in setting II(A) and false in setting II(B), and the assumptions made by PS estimator II are false in both settings. As in settings I(A) and I(B), PS estimator I is approximately unbiased when its parametric assumptions hold [setting II(A)] but is considerably biased for some estimands when its parametric assumptions do not hold [setting II(B)]. For setting II(A), PS estimator II provides slightly biased but low variance estimates of the average effects of treatment across the immune and treatment protective strata and the average effect of treatment across the treatment harmful and doomed strata (the value listed in the Truth column for the PS II estimator for immune=treatment protective is the average effect of treatment among subjects with $S^1 = 0$ and for treatment harmful=doomed is the average effect of treatment among subjects with $S^1 = 1$). For setting II(B), PS estimator II's estimate of the average effect of treatment across the treatment harmful and doomed strata is substantially biased.

In summary, all of the estimators are somewhat sensitive to their assumptions. The PS I estimator makes no assumption about how the average effects within a principal stratum $[E(Y|A=1,X,\underline{S}) - E(Y|A=0,X,\underline{S})]$ compare between principal strata, but is sensitive to its parametric assumptions. The PS II estimator and observed auxiliary stratification

TABLE 3
*Parameter values for simulation study*

| Setting | I(A) | I(B) | II(A) | II(B) |
|---|---|---|---|---|
| $E(Y|A=1, X=0, \underline{S} = immune)$ | 2 | 2 | 2 | 2 |
| $E(Y|A=0, X=0, \underline{S} = immune)$ | 1 | 1 | 1 | 1 |
| $E(Y|A=1, X=0, \underline{S} = treatment\ protective)$ | 2.5 | 2.5 | 2.5 | 2.5 |
| $E(Y|A=0, X=0, \underline{S} = treatment\ protective)$ | 1.5 | 1.5 | 1 | 1 |
| $E(Y|A=1, X=0, \underline{S} = treatment\ harmful)$ | 1.25 | 1.25 | 1.25 | 1.25 |
| $E(Y|A=0, X=0, \underline{S} = treatment\ harmful)$ | 0.75 | 0.75 | 0.75 | 0.75 |
| $E(Y|A=1, X=0, \underline{S} = doomed)$ | 1.75 | 1.75 | 2.25 | 2.25 |
| $E(Y|A=0, X=0, \underline{S} = doomed)$ | 1.25 | 1.25 | 1.25 | 1.25 |
| Distribution within each $A, X, S^0, S^1$ stratum | Normal | Gamma | Normal | Gamma |

TABLE 4
*Simulation study results*

*Settings* I(A) *and* I(B)

| Estimator | Setting I(A) | | | Setting I(B) | | |
|---|---|---|---|---|---|---|
| | Truth | Mean | SD | Truth | Mean | SD |
| Principal stratification estimator I average effects of treatment | | | | | | |
| Immune | 1 | 0.97 | 0.23 | 1 | 0.11 | 0.14 |
| Treatment protective | 1 | 1.00 | 0.13 | 1 | 0.82 | 0.05 |
| Treatment harmful | 0.5 | 0.54 | 0.48 | 0.5 | 0.68 | 0.06 |
| Doomed | 0.5 | 0.51 | 0.19 | 0.5 | 1.36 | 0.10 |
| Principal stratification estimator II average effects of treatment | | | | | | |
| Immune=treatment protective | 1 | 1.00 | 0.03 | 1 | 0.78 | 0.14 |
| Treatment harmful=doomed | 0.5 | 0.50 | 0.04 | 0.5 | 0.60 | 0.08 |
| Observed auxiliary stratification | | | | | | |
| $\Psi_0$ | 1 | 1.00 | 0.12 | 1 | 1.00 | 0.12 |
| $\Psi_1$ | 0.5 | 0.50 | 0.21 | 0.5 | 0.50 | 0.21 |

*Settings* II(A) *and* II(B)

| Estimator | Setting II(A) | | | Setting II(B) | | |
|---|---|---|---|---|---|---|
| | Truth | Mean | SD | Truth | Mean | SD |
| Principal stratification estimator I average effects of treatment | | | | | | |
| Immune | 1 | 0.98 | 0.22 | 1 | 0.00 | 0.14 |
| Treatment protective | 1.5 | 1.52 | 0.14 | 1.5 | 1.91 | 0.06 |
| Treatment harmful | 0.5 | 0.60 | 0.54 | 0.5 | 1.16 | 0.08 |
| Doomed | 1 | 0.96 | 0.21 | 1 | 0.43 | 0.08 |
| Principal stratification estimator II average effects of treatment | | | | | | |
| Immune=treatment protective | 1.42 | 1.33 | 0.03 | 1.42 | 1.63 | 0.04 |
| Treatment harmful=doomed | 0.93 | 0.88 | 0.04 | 0.93 | −0.01 | 0.07 |
| Observed auxiliary stratification | | | | | | |
| $\Psi_0$ | 1.42 | 1.51 | 0.12 | 1.42 | 1.52 | 0.13 |
| $\Psi_1$ | 0.93 | 0.56 | 0.23 | 0.93 | 0.53 | 0.23 |

estimator both put the same constraints on how the average effects within principal strata compare, and both estimators are sensitive to these constraints holding. The PS II estimator was considerably more efficient than the observed auxiliary stratification estimator when its parametric assumptions held, but showed considerable bias for some estimands when its parametric assumptions did not hold. We also performed the same set of simulations for a sample size of 500; the results are not shown but are available from the authors. The performance of the estimators and the way in which the estimators compare were similar to what is shown in Table 3 for the larger sample size of 5000. Some notable features of the simulations for the smaller sample size of 500 were: (1) in a very small proportion of simulations (approximately 0.005%), the observed auxiliary stratification estimator produced estimates that were very large in absolute value (more than 100 times the true value); (2) in a very small proportion of simulations (approximately 0.1%), the PS I estimator did not converge; and (3) the PS estimators exhibited a small amount of bias (a maximum of about 0.1) for some estimands even when their parametric assumptions held.

For the settings considered in Table 3 (with a sample size of 5000), we also computed the expected auxiliary stratification estimator based on model (3). The means of the estimates of $\beta_A + \beta_{\mu A}$ were −3.00, −2.99, −2.18 and −2.19 for settings I(A), I(B), II(A) and II(B), respectively. These mean estimates are not close to the effect of treatment for the subgroup of subjects with $S^1 = 1$. These results illustrate the point mentioned at the end of Section 4.2 that the expected auxiliary stratification estimates cannot necessarily be extrapolated to reliably estimate the effect of treatment for the subgroup of subjects who will have a certain value of the auxiliary variable if treated.

## 6. ANALYSIS OF THE HIP DATA

The Health Insurance Plan (HIP) Study (Shapiro et al., 1988) was a randomized trial of screening for breast cancer, in which more than 60,000 women aged 40–64 at the start were randomized into two groups. In the treatment arm, women were offered screening examinations, consisting of mammography and clinical exam, to be provided in an initial visit and three annual follow-up visits. About one-third of women in the treatment arm refused screening; some of the others did not receive some of the follow-up exams; this information is also recorded. Women in the control group received usual care. The study recorded information on date of death for women who died; cause of death was classified as being due to breast cancer or not. The study also recorded information on whether and when a participant was diagnosed with breast cancer, and if so, whether the cancer was detected by one of the screening exams in the study.

We analyzed the HIP data using G-estimation of structural models. For these analyses, our outcome $Y$ is the natural logarithm of our failure-time, death from breast cancer. Because most subjects in the study do not experience the outcome of interest before the end of follow-up, the outcome is said to be censored.

We consider several choices of the auxiliary outcome $S$. For the first, $S = 1$ if a subject receives a screen, 0 otherwise; for the second, $S = 1$ if a subject is screened and diagnosed with breast cancer, 0 otherwise; for the third, $S = 1$ if a subject is diagnosed with breast cancer and that cancer was detected by the screen, 0 otherwise. Here, the treatment of interest $A$ is (randomized) assignment to the screening arm. We assume that assignment to the screening arm has no effect on the outcome unless a subject actually is screened and is diagnosed (for the third choice, diagnosed based on the screen), and so $\Psi_0 = 0$. Because of censoring, we use as our function

$g\{Y^0(\Psi), X\} = \Delta(\Psi)$, where $\Delta(\Psi) \equiv I\{Y(\Psi) < C(\Psi)\}$, $C(\Psi) = \min\{C, C \exp(\Psi)\}$ (for binary $A$, $S$, and scalar $\Psi$), and $C$ denotes a subject's potential censoring time (here ten years). Other references (Joffe, 2001; Mark and Robins, 1993b; Robins, 1992) discuss the more general approach to dealing with administrative or generalized type I censoring. Joffe (1994) provides fuller justification for use of this particular function. To deal with competing risks, we weight the estimating functions for subjects not censored by competing risks by the inverse of their probability of being uncensored at the end of their follow-up; Robins et al. (1992) provide more details. We consider subjects who died of other causes as censored by a competing risk, although this is somewhat problematic (Kalbfleisch and Prentice, 2002).

For our analysis, we consider the first ten years of follow-up for each subject. In this period, 1259 subjects were diagnosed with breast cancer (626 of these were in the screening arm). There were 340 deaths attributed to breast cancer: 193 in the control arm and 147 in the screening arm. Among the 441 cancers diagnosed among subjects who received at least one screen, 132 cancers were detected by the screen; among the 100 screened breast cancer cases who died of breast cancer, 27 had been detected by the screen.

We applied the modified G-estimation approach to these data. Figure 2 plots the $Z$-test statistics for each choice of $S$ against $\Psi_1$. The test statistics for $\Psi_1 = 0$ are identical $[Z = -2.39, p = 0.017$ (2-sided)]. When $S$ denotes diagnosis of breast cancer among screened subjects, the absolute value of the statistic is minimized around $\Psi_1 = -0.20$, and the 95% confidence interval is $(-0.49, -0.06)$; this means that screening lengthens the time to breast cancer mortality by a factor of $[\exp\{-(-0.2)\} - 1] * 100\% = 22\%$ (6%, 63%) among screened subjects diagnosed with breast cancer. We obtain identical results when $S$ denotes receiving a screen, and so the model estimates may also be interpreted as the effect of screening on screened subjects. This is so because, for any given value of $\Psi$, $g\{Y^0(\Psi), X\} = \Delta(\Psi)$ is the same whether $S$ denotes screening or cancer diagnosis following screening. Because screening is unavailable in the control arm, these estimates are also estimates of complier average causal effects (Imbens and Rubin, 1997). When $S$ denotes breast cancer diagnosis as the result of a screen, the corresponding point estimate (95% confidence interval) is

−0.76 (−1.15, −0.27); the corresponding lengthening is 114% (31%, 216%). Thus, the effect of screening among screened subjects diagnosed with breast cancer is less than the effect on the subset of these subjects whose diagnosis was a result of the screening, as would be expected.

## 7. UTILITY OF VARIOUS APPROACHES

Causal modeling has multiple purposes. These include assisting in making decisions about possible interventions, predicting the results of those interventions, and better understanding of the processes leading to the outcome(s) under study. We consider the approaches sketched above in terms of their utility for these purposes, as well as generalization of study results.

### 7.1 Making Decisions

Making decisions about possible interventions and predicting the results of those interventions are closely linked. Normally, one would want to choose, for any individual, the treatment that leads to the best possible expected result, however that is defined. This decision must be made on the basis of information available at the time of the decision; for point exposures, that means only baseline covariates may be used to predict outcomes under a given treatment and so guide treatment decisions.

We formalize this within a decision-theoretic framework. Let $L_i(S^a, Y^a; a)$ denote an individual's loss function associated with decision $a$. The loss function may be associated with the subject's principal stratum, observed auxiliaries or outcomes, measured pretreatment covariates, or other individual-specific factors. To formalize this idea, let $Q$ denote a collection of variables. We can write the expected loss given $Q$ as $E\{L_i(S_i^a, Y_i^a, a)|Q\}$. Under the usual decision-theoretic framework, a decision should depend only on information available at the time of the decision; further, if the decision is to be based only on information available in the statistical models fit to the data, we can condition only on measured baseline covariates $X$. Under these constraints, decisions should be based only on comparisons of $E\{L_i(S_i^a, Y_i^a, a)|X\}$ for different treatments $a$, and so depend on the joint density of the potential outcomes $f(\underline{S}, \underline{Y}|X)$ only through the marginal densities $f(S^a, Y^a|X) = f(S^a|X)f(Y^a|S^a, X)$, where $\underline{S} \equiv \{S^a\}$, $a \in \mathbf{A}$ and $\mathbf{A}$ denotes the set of possible treatments $a$. Thus, the optimal feasible decision based solely on data available to the study will depend on the marginal distributions $f(S^a|X)$ and $f(Y^a|S^a, X)$, both of which are estimable using conventional methods (the latter as described in Section 2.2); the use of the factorization of the marginal density $f(S^a, Y^a)$ for estimation will be particularly important when $Y^a$ is undefined for some values of $S^a$ (e.g., with censoring by death). Relying solely on the marginal density $f(Y^a|X)$ will be adequate for deciding among treatments if the loss does not depend on $S^a$ or if the auxiliary variable is not affected by treatment and the loss function may be written as the sum of separate components associated with main and auxiliary outcomes [i.e., if $L_i(S^a, Y^a; a) = q_1(S^a) + q_2(Y^a)$ for some functions $q_1(\cdot)$ and $q_2(\cdot)$]. If the loss is a function of $S^a$, none of the estimands discussed in Section 3 is adequate by itself, as they focus only on the density of $Y^a$ [some methods proposed in conjunction with PS (Imbens and Rubin, 1997) yield estimates of $f(S^a|X)$ and so potentially could be used to estimate $f(S^a, Y^a)$]. In the nephrology example, the auxiliary outcome, ESRD, is a condition with severe and life-disrupting consequences (requiring dialysis or transplant) which would normally influence the loss function.

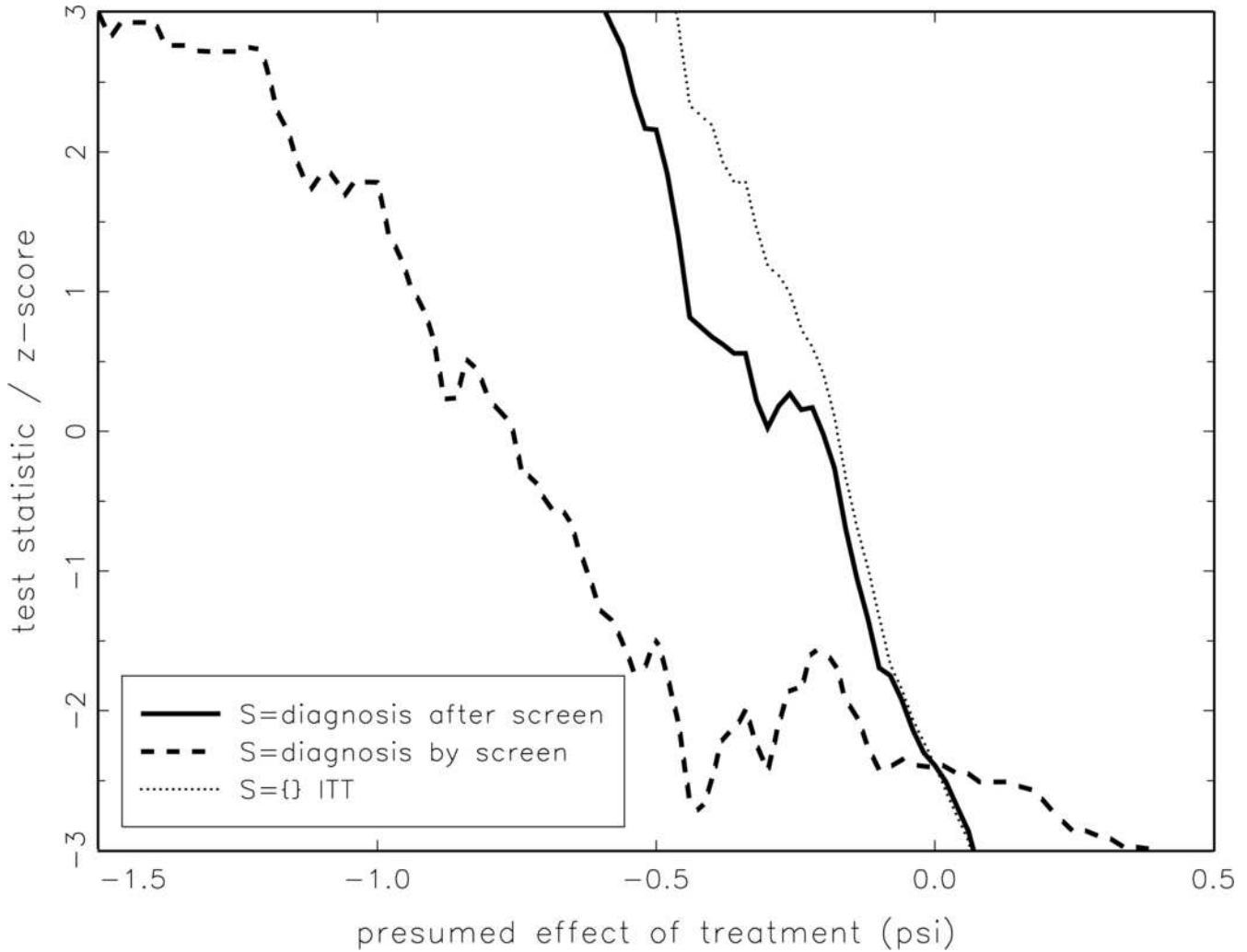

Suppose that the loss is not a function of $S^a$. In the cancer screening example, it is arguable that the benefits of screening should be assessed solely in terms of mortality $Y^a$, even though there are costs that are associated with cancer diagnosis, including the financial costs of administering screening exams, patient inconvenience, and the costs of falsely diagnosing a subject as having cancer who actually does not, which can both lead to unnecessary surgical procedures and have adverse psycho-logical impact (Brett, Bankhead, Henderson, Watson and Austoker, 2005; Barratt, Irwig, Glasziou, Salkeld and Houssami, 2002). If the loss is a function of $Y^a$ alone, optimal decisions should be based on $f(Y^a|X)$; that is, we should seek to minimize $E\{L_i(S^a, Y^a, a)|X\} = E\{L_i(Y^a, a)|X\} = \int_{Y^a} L_i(Y^a, a)f(Y^a|X)\,dY^a$. Suppose that, as is typical, $L_i(Y_i^a, a)$ is a monotonic function of $Y_i^a$ for all subjects (e.g., mortality is never a desired outcome); suppose further that $L_i(Y_i^a, a)$ does not depend on $a$ and so $L_i(Y_i^a, a) = L_i(Y_i^{a'}, a')$ for all $a, a'$ (e.g., the loss associated with mortality does not depend on whether one had been screened). Then, the relative expected loss for an individual under one treatment relative to another can be assessed by comparing the densities $f(Y^a|X)$ and $f(Y^{a'}|X)$.

Although these quantities may be evaluated from a model which conditions on a post-treatment auxiliary, such evaluation will typically be more involved than the simpler and more direct modeling of the marginal densities $f(Y^a|X)$. Let $Z$ denote a post-treatment auxiliary; $Z$ can be a principal stratum or the single potential auxiliary $S^a$. We can compute and estimate $f(Y^a|X, Z)$ using methods discussed above (including the methods based on observed auxiliary stratification; Section 4.3); nonetheless, to make a decision, we must integrate out the auxiliary $Z$, as it is unknown at the time a decision is being made. Thus, from models based on PS or single potential auxiliaries, we can assess the loss by evaluating $f(Y^a|X) = \sum_z f(Y^a|X, Z = z)pr(Z = z|X)$; these calculations are more involved than those based on stratifications based on the observed covariates $X$, where the most relevant comparison of densities can be read simply (or directly) from a single regression of $Y$ on $X$ and $A$. Further, the increased number and complexity of the models required to evaluate the desired quantities may lead to increased likelihood of obtaining incorrect conclusions due to model misspecification. If the sharp null hypothesis of no effect is true, tests of the null using G-estimation of parameters in structural nested mean models will be (relatively) robust to model misspecification in randomized trials, because they essentially are based on intent-to-treatment tests of the null (Robins, 1992).

Information on the expected value of post-treatment auxiliaries more generally might be of use in making decisions in situations in which decision makers had a better idea of the value of these post-treatment auxiliaries than can be assessed from the data. This could happen if, for example, a new test were developed that discriminates well between subjects who will and will not develop the auxiliary outcome. Then, effects for groups defined by an auxiliary variable might be used to infer the effects of treatment for groups defined by this baseline variable. Additionally, finding that effects of treatment differ for subjects with different levels of the auxiliary variable could be a useful spur to find baseline variables which predict the auxiliary variable well.

For making decisions, even the expected auxiliary outcomes play no special role in principle. For this purpose, covariates which strongly modify the effects of treatment (measured on a clinically relevant scale) are of greatest value. Sometimes, the expected auxiliary may be such an effect modifier; in the cancer screening example, it is plausible that the risk of breast cancer would be an important modifier of the benefit of screening, especially if the benefit is assessed on the scale of difference in probability of mortality if screened or not. Unfortunately, we cannot evaluate this in the HIP study, which collected little covariate information.

7.1.1 *Treatments given over an extended period.* Where the treatment is actually applied over an extended period of time, estimands stratifying on post-treatment auxiliaries may have more immediate relevance for decision making. Suppose that the auxiliary variable is measured shortly after the initiation of treatment. Suppose further that treatment never changes over the course of follow-up and that the effect of treatment on the ultimate outcome, measured at the end of a fixed follow-up period, is cumulative (i.e., the effect of treatment given early during follow-up on this outcome is in the same direction as the effect of treatment later). Then, finding that treatment is beneficial for some groups defined by post-treatment variables and harmful for others will lead to recommendations to discontinue treatment for the group for whom treatment is harmful.

To illustrate and formalize this, consider elaborating model (4) to model the effect of the different components $A_k$ of treatment received at different times $k, k = 1, \ldots, K$. Let $\dot{A} \equiv (\sum_{k=1}^{K} A_k)/K$ be the average value of treatment received during follow-up. The model

$$(7) \qquad Y^0 = Y - \dot{A}(1 - S)\boldsymbol{\Psi}_0 - \dot{A}S\boldsymbol{\Psi}_1$$

is indistinguishable from model (4) when treatment for all subjects remains constant over time; if the model is correct and $S$ is assessed early in follow-up, one would then be justified in making treatment decisions beyond $k$ based on the value of $S$ measured by $k$. Such inference is, however, strongly dependent on modeling assumptions. To see this, consider the model

$$(8) \qquad Y^0 = Y - A_1(1 - S)\boldsymbol{\Psi}_0 - A_1 S \boldsymbol{\Psi}_1,$$

which specifies that the only treatment which affects the outcome is that received during the first period. When treatment $A_k$ remains constant over time, this model is indistinguishable in the data from model (7); nonetheless, one model will suggest that treatment may be beneficial beyond $k$, whereas the other will not.

## 7.2 Explanatory Analyses

Another possible use of analysis stratifying on post-treatment auxiliaries is explanatory. In scientific work, one typically wants to explain the data in ways which will lead not only to an understanding of the data themselves but also to generalization to other settings (often beyond the sampling frame of the study). In this section, we consider the utility of analysis stratifying on post-treatment auxiliaries for such explanation and contrast it with approaches which concentrate on causal mechanism.

7.2.1 *Two types of explanation: effect modification vs. causal mechanism.* Effects which vary by post-treatment auxiliaries can be used to explain observable differences in outcomes between groups receiving different treatments (White and Goetghebeur, 1998). For example, in the HIP study, an effect of screening $\boldsymbol{\Psi}_1$ of $-0.2$ for screened subjects diagnosed with breast cancer (or $-0.76$ for screened subjects diagnosed due to the screen), together with no screening effect in the remaining subjects, can explain the differences in survival between the randomized groups. We say that there is effect measure modification (Rothman and Greenland, 1998) here by the post-treatment auxiliary (diagnosis after screen or diagnosis due to screen) because the measure of treatment effect differs between the different groups defined by the auxiliary [i.e., $\boldsymbol{\Psi}_1 \neq \boldsymbol{\Psi}_0$ in (4)]; this variation in effects can be used to explain the differences between the groups.

One might be tempted to conclude that cancer diagnosis (or its correlates) participates in the processes leading to improved survival and perhaps interacts causally with the treatment of interest. Such mechanistic inference will often be suggestive but is fraught with pitfalls (Thompson, 1991). In particular, the concept of effect measure modification (a.k.a. statistical interaction) differs from mechanistic interaction in that effect measure modification does not require that the modifier itself be a variable which can be modified directly, whereas mechanistic interaction requires the modifier to be the subject of intervention. Thus, one need not speak of the causal effect of an effect modifier.

Effect modification by post-treatment auxiliary variables provides less satisfactory and less satisfying explanations of observed associations when divorced from the concepts of causal mechanism and causal intermediates. To see this, consider a scenario in which the effect of aggressive blood pressure management on MI is greater for subjects who subsequently develop ESRD [e.g., $|g\{E(Y^1|A=1,S=1)\} - g\{E(Y^0|A=1,S=1)\}| > |g\{E(Y^1|A=1,S=0)\} - g\{E(Y^0|A=1,S=0)\}|$ for some monotone link function $g(\cdot)$; typically g(y) is either $y$ (the identity link, leading to risk differences), ln($y$) (the log link, leading to risk ratios) or logit(y) (the logit link, leading to odds ratios)]. Contrast the following four statements:

1. The effect of aggressive treatment of blood pressure on MI is greater for people who subsequently develop ESRD (observed auxiliary stratification), with the effects in the different subgroups explaining the overall difference between subjects treated aggressively and those not.
2. The effect of aggressive treatment of blood pressure on MI is greatest for people who would develop ESRD only if not treated aggressively (principal stratification).
3. Aggressive treatment of blood pressure prevented the development of ESRD in some subjects (one principal stratum). This, in turn, prevented the development of MI for some subjects. Thus, ESRD mediated in part the effect of aggressive management of blood pressure on MI.
4. Aggressive treatment of blood pressure prevents the development of ESRD in some subjects. This, in turn, prevents the development of MI for some subjects. Thus, ESRD mediates in part the effect of aggressive management of blood pressure on MI.

The last two statements attempt to provide some understanding of the causal mechanisms leading from blood pressure treatment to MI; these are often considered in terms of direct and indirect effects (Pearl, 2001; Robins and Greenland, 1992; Robins, 1999). The first two do not attempt such explanation. Statements 3 and 4 are identical except for the tense. Statement 3 is an attempt to explain what has happened; statement 4 refers to more general causal processes and in principle is an attempt to predict the future and so is more ambitious.

In general, we prefer statements 3 and 4. Statements 1 and 2 are not, in general, useful for making treatment decisions, and the quality of their explanations of the data is poor. Where warranted based on appropriate subject-matter considerations, statements 3 and 4 provide informative explanations. At the risk of generalizing beyond the data at hand, statement 4 makes general statements about causal processes which are then potentially testable in other data. An important role of science is to extrapolate from one's data and make predictions which may be testable in other data or designs; this is most easily fostered by considering causal mechanisms. This important scientific goal is sometimes fostered by considering effects of auxiliary variables that are not under the direct control of the investigator in the given study; Rubin (2004) presents a somewhat contrary view.

In the nephrology example, it is meaningful to consider the effect of ESRD on MI; this accords nicely with common usage. The effect can be approximated by comparing what would happen to someone who has or develops ESRD with what would have happened had that person received a transplant from a genetically identical person (an identical twin or identical triplets) who did not have ESRD. Even though this experiment could rarely, if ever, be carried out, it provides a useful approach for defining the effect in terms of a thought experiment.

DAGs provide a nice intuitive representation of direct and indirect effects; the ideas of potential outcomes and counterfactuals allow these ideas to be made more precise. Figure 3 shows causal relationships among the variables. In this graph, the path $A \rightarrow S \rightarrow Y$ represents the indirect effect of $A$ on $Y$ (i.e., that part of the effect of $A$ that is mediated by the specified variable $S$; it is necessary to specify the auxiliary variable(s) $S$ to define what is meant by indirect and direct effects), and $A \rightarrow Y$ represents the direct effect of $A$ on $Y$ (i.e., that part of the effect not mediated by $S$). Graphs like this are sometimes known as path diagrams and have been used to justify the use of linear models for multivariate normal data (Pearl, 2000). The path-analytic approach suffers from the lack of a nonparametric definition of causal effects and generally unrealistic assumptions of multivariate normality; further, it does not extend naturally to settings with interactions among the variables or to nonlinear models.

The ideas of potential outcomes may be used to define direct, indirect and joint effects of treatment. Let $Y^{a,s}$ denote the (continuous) outcome one would see for a given subject at level $a$ of the main treatment of interest (e.g., management of hypertension) and level $s$ of the auxiliary outcome (e.g., ESRD). Underlying the notation is the assumption that the auxiliary outcome is in some sense manipulable. The notation has been used in conjunction with PS (Angrist, Imbens and Rubin, 1996) as well as discussions of direct and indirect effects.

The direct effect of the main treatment $A$ controlling for an auxiliary variable $S$ may be defined as the contrast of $Y^{a^1,s}$ and $Y^{a^0,s}$ for an individual or a group; that is, the direct effect contrasts the outcome that would be seen under two different levels of the primary treatment, physically manipulating the auxiliary variable to a given level $s$. Direct effects are not uniquely defined; there are separate direct effects for each level of the auxiliary variable $s$.

There are, in fact, two separate types of direct effects (Pearl, 2001; Robins and Greenland, 1992; Robins, 2003), depending on the nature of the manipulation of the auxiliary variable. If the auxiliary variable is set to a prespecified level $s$, the resulting contrasts are known as "controlled" direct effects. If the auxiliary variable is set to the level it would have reached at some reference level $a^*$ of the treatment, the effects are known as "natural" direct effects; natural direct effects are contrasts of $Y^{a^1,S^{a^*}}$ and $Y^{a^0,S^{a^*}}$, $a^* \in \{a^0, a^1\}$. For $a^* = a^0$, the natural direct effect is a contrast of the potential outcome that would have been seen had the main treatment been set to level $a^1$ but the auxiliary variable set to the level $S^{a^0}$ it would have taken had the subject been given level $a^0$ with the potential outcome that would have been seen had the subject been given level $a^0$ of the main treatment.

Although there is no general definition of controlled indirect effects, natural indirect effects are nonparametrically defined as contrasts of $Y^{a^*,S^{a^1}}$ and $Y^{a^*,S^{a^0}}$; that is, for $a^* = a^0$, the indirect effect is a contrast of what would have happened had a subject been given treatment $a^0$ but had his or her auxiliary variable set to the level $S^{a^1}$ it would have attained had he or she been given $a^1$ with what would have happened had the subject been given $a^0$. Statements 3 and 4 above may be understood as speaking about the natural indirect effects of aggressive management of blood pressure.

### 7.2.2 *Contrasting models for contrasting explanations.*

We present here a simple semiparametric model which allows us to contrast the ideas of mediation with those presented previously, in particular, with the modification of treatment effect by post-treatment covariates, as in (4). We consider a continuous outcome variable $Y$. A simple model for the joint effects of hypertension management and ESRD is

$$(9) \qquad Y^{a,s} = Y^{0,0} + a\gamma_1 + s\gamma_2 + as\gamma_3.$$

In this model, $\gamma_1$ represents the controlled direct effect of aggressive treatment of blood pressure (holding $S$ at 0), $\gamma_2$ represents the effect of modification of ESRD and $\gamma_3$ represents an interaction between ESRD and aggressive treatment of blood pressure. Under an assumption of monotonicity of the effect of treatment $A$ on $S$ (e.g., $S^1 \le S^0$), the parameter $\gamma \equiv \{\gamma_1, \gamma_2, \gamma_3\}$ is related to the effects of aggressive treatment of blood pressure on principal strata as follows:

for subjects doomed to develop ESRD ($S^1 = S^0 = 1$), the effect of aggressive treatment ($Y^1 - Y^0 = Y^{1,S^1} - Y^{0,S^0}$) is $\gamma_1 + \gamma_3$;

for subjects immune to ESRD ($S^1 = S^0 = 0$), the effect of aggressive treatment is $\gamma_1$;

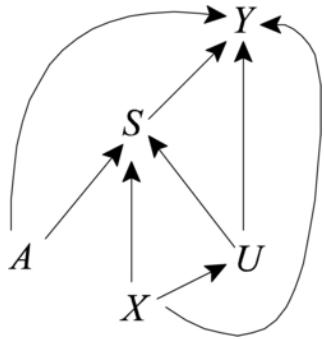

Fig. 3. *A directed acyclic graph representing the relations among the variables in the example of Section* 7.2.

for subjects for whom treatment prevents ESRD, the effect of aggressive treatment is $\gamma_1 - \gamma_2$.

In this model, statement 2 (now relating to a continuous cardiovascular outcome) that the beneficial effect of aggressive treatment of blood pressure is greatest for people who would develop ESRD only if not treated aggressively is implied by $\gamma_2 > 0$, $\gamma_3 > -\gamma_2$ (assuming small values of $Y$ are better outcomes). Further, if $\gamma_2 > 0$ and aggressive treatment and ESRD do not interact in affecting the outcome (i.e., $\gamma_3 = 0$), then statement 2 is true, as is statement 4.

Further, the effect of treatment among treated subjects who do not develop ESRD is

$$E(Y^1|A = 1, S = 0) - E(Y^0|A = 1, S = 0)$$
$$= E(Y^{1,S^1} - Y^{0,S^0}|A = 1, S^1 = 0)$$
$$= \sum_s pr(S^0 = s|A = 1, S^1 = 0)$$
$$\quad \cdot E(Y^{1,S^1} - Y^{0,S^0}|A = 1, S^1 = 0, S^0 = s)$$
$$= \gamma_1 - pr(S^0 = 1|S^1 = 0)\gamma_2,$$

and the effect of treatment among the treated who develop ESRD is $\gamma_1 + \gamma_3$. Thus, $\mathbf{\Psi}_0$ in (4) equals $\gamma_1 - pr(S^0 = 1|S^1 = 0)\gamma_2$, and $\mathbf{\Psi}_1 = \gamma_1 + \gamma_3$.

In a randomized trial of $A$, the parameters in model (9) are identified from the data under two types of assumptions:

1. The assumption that the initial treatment $A$ is randomized, along with the modeling assumptions inherent in (9) (Robins and Greenland, 1994; Ten Have et al., 2004). In this approach, the modeling assumptions are not fully testable, even in very large studies. This approach is similar to that proposed above for observed auxiliary stratification and under some conditions is essentially identical. For G-estimation, the approach requires the same number of estimating equations as the number of free parameters. Under the above normality and homoscedasticity assumptions and the assumption that $f(Y^{0,0}|X, A, S^1) = f(Y^{0,0}|X)$, the optimal function is the vector function $g\{Y^{0,0}(\gamma), X\} = \varepsilon(\gamma)\{1, E(S|X, A = 1) - E(S|X, A = 0), E(S|X, A = 1)\}^T$. We now require the covariate vector to predict not only the level of the auxiliary variable among the treated [as was true for model (4)], but also to predict the effect of the primary treatment on the auxiliary, and we require the three elements $\{1, E(S|X, A = 1) - E(S|X, A = 0), E(S|X, A = 1)\}$ to not be collinear; this often requires at least two covariates. Restrictions on the parameter space can result in a smaller vector of estimating functions; typically, the restrictions are either that there is no interaction between $A$ and $S$ (i.e., $\gamma_3 = 0$) (Robins and Greenland, 1994; Ten Have et al., 2004), or that some parameters play no role in explaining the data and so are not estimable. To explain the latter, suppose that one level of a binary auxiliary variable $S$ does not occur for some level of treatment $A$; for example, $pr(S = 1|A = 0) = 0$; in the HIP study, this is true both when $S$ denotes screening and when $S$ denotes diagnosis due to screen. Under this condition, $S = AS$ for all subjects, and one cannot simultaneously estimate both the main effect of $S$ ($\gamma_2$) and its interaction with $A$ ($\gamma_3$) in (9). Under this model, for a binary treatment and auxiliary variable, the effect of the auxiliary treatment $S$ ($\gamma_2$) is the same as the effect of the main treatment $A$ among treated subjects with $S = 1$($\mathbf{\Psi}_1$). Thus, the effect of assignment to screening for screened women who are diagnosed due to the screen [$\mathbf{\Psi}_1$ in (4)] might also be interpreted as the effect of diagnosis or early detection by screen [$\gamma_2$ in (9)], and the effect of assignment to screen among subjects who received the screen might also be interpreted as the effect of screening.

2. The assumptions above, plus the assumption the assignment of the auxiliary variable $S$ is ignorable; that is, $pr(S|X, A, L, Y^{a,s}) = pr(S|X, A, L)$, where $X$ refers to covariates measured at baseline and $L$ to covariates measured after baseline but before $S$ (Robins and Greenland, 1992). Under the additional but untestable ignorability assumption, the structural model (9) is fully testable.

The DAG in Figure 3, which corresponds to a randomized trial of the main treatment $A$, is consistent with assumption 1 but not assumption 2. The assumption of initial randomization is justified because there are no arrows into $A$; further, the arrow from $X$ to $S$ is consistent with an association between $X$ and $S$ among treated subjects (also part of assumption 1). Assumption 2 is not justified, because the arrows from $U$ to $S$ and $U$ to $Y$ imply that the effect of $S$ on $Y$ is confounded.

The joint effects of $A$ and $S$ are, in principle, identifiable without making modeling assumptions in an experiment in which both factors (e.g., blood pressure treatment and kidney function) are experimentally controlled and both treatments assigned ran-domly. Kidney function might be controlled experimentally by surgery. Ethical considerations preclude performing such experiments in humans; animal experimentation could, in principle, be used to learn about the direct effect of blood pressure control.

Model (9) is fairly simple and might not be a faithful representation of the real-world situation. In particular, the effects $A$ and $S$ in (9) might vary with observed covariates (e.g., baseline covariates $X$ or observed treatment levels $A$) or latent covariates (e.g., the principal strata). Approaches for dealing with these issues (Robins, 1999; Robins, Rotnitzky and Scharfstein, 2000; van der Laan and Petersen, 2004; Robins, 2003) are beyond the scope of this paper.

7.2.3 *Explanation and generalization.* Explanatory analysis can serve to explain the findings in the data at hand in ways which may generalize further to other populations or settings, a more ambitious task requiring further assumptions. Generalizability beyond one's data is enhanced if the relations found in one's data are likely to hold in other settings. Two strategies for this are:

1. Obtaining a finer understanding of the causal processes leading from treatment to outcome. With sufficiently detailed and accurate understanding, altering one component in the causal system may not change the other causal mechanisms in the system; this may allow better prediction of the effects of interventions in other settings (Pearl, 1995, 2000). Partitioning effects into component direct and indirect effects is an attempt to obtain understanding in those terms.

2. Estimating the effects of treatment for homogenous subgroups in which causal effects in the population under study can be assumed to be similar to like groups in other populations. This can be done by estimating the effects of treatment for subgroups defined by observable variables $X$ or $S$, or by the latent principal strata $\underline{S}$. Principal stratification can lead to a finer partition of the population than stratification on observed variables alone, and so can potentially lead to more generalizable conclusions.

These strategies are not mutually exclusive. Structural nested models have been formulated in a way which recognizes possible differences in effect in subgroups not identifiable at baseline (Robins et al., 1992; Robins, Rotnitzky and Scharfstein, 2000); our approach to using G-estimation for estimating treatment effects for strata defined by observed post-treatment auxiliaries allows finer stratification. Further, PS approaches have been used in conjunction with approaches which allow one to consider partitioning causal effects (Angrist, Imbens and Rubin, 1996).

## 8. FUTURE WORK: EXTENSIONS OF ESTIMANDS AND ESTIMATION

We conclude the paper with a discussion of various additional complications which may arise in applied problems, and outline possible directions for future work in this context.

### 8.1 Multiple Types of Auxiliary Variables

In considering conditioning on a post-treatment auxiliary $S$, we have so far considered settings in which the outcome of interest $Y$ is meaningful for all levels of $S$. Sometimes, the outcome of interest $Y$ is not defined at some levels of $S$. Most notably, if $S$ is (or includes) an indicator of vital status, the outcome $Y$ may not be defined meaningfully for subjects who die, and so the effect of $A$ on $Y$ (i.e., the contrast of $Y^1$ and $Y^0$) will be defined only for subjects who would live whether or not treated (i.e., the principal stratum in which $S^1 = S^0 = 1$, where $S = 1$ indicates being alive).

In the nephrology example, some subjects may develop ESRD only if their blood pressure is treated aggressively, because aggressive treatment may prevent mortality. Thus, a simple version of the monotonicity assumption stating that aggressive treatment never causes ESRD will not be true. Because of countervailing biases, the direction of bias of the naive approach (comparing among treated and untreated subjects who develop ESRD) may not be predictable.

Where there is a composite auxiliary variable (e.g., death and ESRD), there are several possible estimands. Let $S_1$ denote an auxiliary variable for which outcomes are well defined at all values of the variable (e.g., ESRD), and let $S_2$ denote an auxiliary variable for which outcomes are well defined only for one value of the variable (e.g., death; $S_2 = 1$ if a subject lives, 0 if dead). Because other outcomes are not meaningfully defined for people who die (Frangakis and Rubin, 2002; Kalbfleisch and Prentice, 2002), it is difficult to ignore mortality in defining causal estimands for this setting. This leaves fewer options for dealing with mortality: looking at the effect of treatment on the conditional distribution of $Y$ given survival, and PS. Within the framework of PS for mortality, almost any of the solutions above for the other auxiliary variable (ESRD) may be applied for defining causal estimands. The possible estimands include comparisons of the following expectations for different levels $a$ of treatment:

1. $E(Y^a|S_2^0 = 1, S_2^1 = 1)$, the effect of treatment for subjects who would not die whether or not treated.
2. $E(Y^a|S_2^a = 1)$, the effect of treatment on the expectation of MI for subjects who would not die under that treatment. As before (Section 2.2), this is not a comparison of outcomes for a common set of subjects. This estimand is most easily understood in conjunction with the effect of treatment on survival [i.e., comparisons of $E(S_2^a)$]. This and the previous estimand are not defined in terms of ESRD. Estimands which use ESRD include:
3. $E(Y^a|S_1^a = 1, S_2^a = 1)$, the effect of treatment on the distribution of MI among subjects who would be alive and have ESRD under that treatment;
4. $E(Y^a|S_1^0 = 1, S_1^1 = 1, S_2^0 = 1, S_2^1 = 1)$, the effect of treatment among subjects who, under both treatment levels, would be alive and develop ESRD (full or dual PS);
5. $E(Y^a|S_1^1 = 1, S_2^0 = 1, S_2^1 = 1)$, the effect of treatment on subjects who would live whether or not treated and develop ESRD if treated (single potential auxiliary stratification for $S^1$ in conjunction with PS for $S^2$);
6. $E(Y^a|S_1 = 1, S_2^0 = 1, S_2^1 = 1, A = 1)$, the effect of treatment on subjects who would live whether or not treated, who are treated and develop ESRD (observed auxiliary stratification); and
7. $E\{Y^a|E(S_1^0|X), S_2^0 = 1, S_2^1 = 1, A = 1\}$, the effect of treatment on subjects at a particular risk of ESRD who would live whether or not treated.

We expect that the likelihood-based or sensitivity analysis methods mentioned in Section 4.4 could be extended to deal with these issues and estimands.

## 8.2 Other Extensions

Many real-world problems, including studying the effect of the aggressive management of blood pressure on renal patients, involve problems that are substantially more complicated than those considered here. Complications arise from the fact that all three main variables under study (i.e., treatment $A$, auxiliary variable $S$ and outcome $Y$) may be more complex than the simple binary variables we have discussed. The additional complexity of each may require refinement or redefinition of the effects under study as well as of the methods used to estimate them.

The auxiliary variable $S$ of interest may, in fact, be measured repeatedly over time. In the renal study, the time ESRD develops will be noted; in the HIP study, the time of breast cancer diagnosis will be noted. For such a failure-time variable $S$, one could define effects of treatment based on the actual failure-time [e.g., compare $E(Y^a|S^1, S^0)$ for different $a$], or on whether the failure-time exceeds some threshold [e.g., compare $E\{Y^a|I(S^1 > s), I(S^0 > s)\}$ for different $a$]. For failure-time outcomes, as in the HIP study, we can revise the causal estimand to account for the timing of changes in the auxiliary variable. Let $S(t) = 1$ after a subject with breast cancer is diagnosed by screen, 0 otherwise (a modification of the third definition of $S$ in Section 6). A version of the accelerated failure-time model with time-varying covariates [a generalization of (4)] is $T^0 = \int_0^T \exp\{AS(t)\Psi\} dt$ (Cox and Oakes, 1984; Robins, 1992; Robins et al., 1992), where $T$ is the subject's failure time and $T^0$ the failure-time if the subject had not been screened. The causal parameter $\Psi$ now represents the effect of screening in shortening time from cancer diagnosis to mortality among screen-diagnosed subjects. G-estimation could be applied for estimating parameters in this model.

Similarly, the outcome variable $Y$ may be a repeated measures outcome. This allows many options for the time that the value of the auxiliary variable is measured: one may be interested in defining effects of treatment conditional on the value of the auxiliary variable at the time the outcome $Y$ is measured, or one year previously, or on the time of failure for a failure-time auxiliary such as ESRD. If the outcome is defined at all levels of the auxiliary, all options are potentially meaningful. We expect that any of the methods for defining and estimating causal quantities sketched above (Sections 3 and 4) could apply.

The study exposure or treatment is often not a simple scalar but may vary over time, and the joint effects of treatments received at different points in time may be of interest; this is true of observational studies of the effect of aggressive management of blood pressure, in which therapy is provided over an extended period and changes may be made over that time. Robins (Robins, Rotnitzky and Scharfstein, 2000) has provided a general approach to defining the component and joint effects of treatments given over an extended period; this approach would need to be generalized to allow the effects of a component of treatment or of a treatment plan to depend on auxiliary variables subsequent to the treatment. For defining the component effect of a treatment $A_t$ applied at $t$, we might want to allow the effect to depend on $S_{t+1}$, for example. Alternatively, we might want to allow the joint effect of the component treatments in a prespecified regime to depend on the level of an auxiliary variable that applies at a fixed point in time after the start of follow-up.

## ACKNOWLEDGMENTS

The authors thank Zhen Chen, Michael Elliott and Tom Ten Have for useful discussions, and the referees and Editors for incisive comments which led to substantial improvements to the manuscript. The authors also thank Tom Riley, Phil Prorok and Brian Irons for providing access to the HIP data. This work was supported in part by National Institutes of Health Grants U01-DK60990 and R01-CA095415.